\documentclass[a4paper]{jpconf}
\usepackage{graphicx}
\usepackage{amsmath,amsthm,amssymb,amsopn,amsfonts}
\usepackage{mathtools}
\usepackage{multirow}
\usepackage{algorithmic}
\usepackage{algorithm}
\usepackage{siunitx}

\newcommand{\degree}{^{\circ}}

\begin{document}

\bibliographystyle{iopart-num}

\title{A Bayesian Approach to CT Reconstruction with Uncertain Geometry}

\author{Frederik H Pedersen, Jakob S J\o{}rgensen, and Martin S Andersen}

\address{Department of Applied Mathematics and Computer Science, Technical University of Denmark, Kongens Lyngby, Denmark}

\ead{\{fhape,jakj,mskan\}@dtu.dk}

\begin{abstract}
Computed tomography is a method for synthesizing volumetric or cross-sectional images of an object from a collection of projections. Popular reconstruction methods for computed tomography are based on idealized models and assumptions that may not be valid in practice. One such assumption is that the exact projection geometry is known. The projection geometry describes the relative location of the radiation source, object, and detector for each projection. However, in practice, the geometric parameters used to describe the position and orientation of the radiation source, object, and detector are estimated quantities with uncertainty. A failure to accurately estimate the geometry may lead to reconstructions with severe misalignment artifacts that significantly decrease their scientific or diagnostic value. We propose a novel reconstruction method that jointly estimates the reconstruction and the projection geometry. The reconstruction method is based on a Bayesian approach that yields a point estimate for the reconstruction and geometric parameters and, in addition, provides valuable information regarding their uncertainty. This is achieved by approximately sampling from the joint posterior distribution of the reconstruction and projection geometry using a hierarchical Gibbs sampler. Using real tomographic data, we demonstrate that the proposed reconstruction method significantly reduces misalignment artifacts. Compared with two commonly used alignment methods, our proposed method achieves comparable or better results under challenging conditions.
\end{abstract}

\section{Introduction}
Computed tomography (CT) is a computational imaging technique for synthesizing volumetric or cross-sectional images of an object from a collection of projections. The projections are obtained by illuminating the object with radiation at different orientations and measuring the attenuation of the radiation due to absorption. The collection of all projections is known as the sinogram, and the process of reconstructing an image of the interior absorption of an object from the sinogram is an ill-posed linear inverse problem. Numerous reconstruction methods for CT have been developed, see e.g. \cite{PCHbook} for an overview of common reconstruction methods including projection methods, algebraic methods, and variational methods. Computed tomography is used in many different scientific fields, the most well known being medical imaging, where it provides a non-invasive method for inspecting the interior of patients \cite{medicalCT}. It is also popular in geoscience, where it is an invaluable tool to analyze micro-structure of materials \cite{GeophysicsCT,MaterialCT}.
\newline
\newline
Experimental CT setups are comprised of a radiation source, an object to be scanned, and a detector measuring the radiation after attenuation by the object. The projection geometry is the relative position of the source, object, and detector. Popular reconstruction methods rely on the assumption that the projection geometry is known. However in practice the projection geometry is estimated and therefore uncertain. Mismatches between the actual projection geometry from the data acquisition, and the estimated projection geometry used in the reconstruction model may lead to reconstructions with severe misalignment artifacts that degrade their quality, and hence also their scientific or diagnostic value.
\newline
\newline
To mitigate misalignment artifacts, correction of the projection geometry for CT is essential. Many different methods have been proposed to address various errors in the projection geometry. One class of methods is based around doing a calibration scan with fiducial markers, and then correcting the geometry based on the acquired scan \cite{CalibrationScan}. However taking additional measurements may be impossible or very costly depending on the specific application or experimental setup. For this reason it is of great interest to have marker-free methods, that do not rely on additional measurements. One class of such methods corrects the center-of-rotation based solely on the measured sinogram without a need to compute additional reconstructions. These methods typically rely on inherent symmetry features in the recorded sinogram. For instance, \cite{CenterOfMass} corrects the center-of-rotation based on a center-of-mass calculation, and \cite{CrossCor,FourierCrossCor} identifies the center-of-rotation by computing cross-correlation between projections. The main advantage of these methods is their simplicity and low computational cost. However, they are highly dependent on the quality of the data, and thus they may suffer if provided with highly noisy data. Another class of methods is based on image metrics in the reconstruction space \cite{TVMetric,ImageMetric}. The correction is done based on computed reconstructions, which makes it more robust against noise (due to possibility of regularization). The downside of these reconstruction-based methods is increased computational cost, and they may also suffer if the reconstruction contain other artifacts, e.g., ring-artifacts due to uncertainty in the source intensity. Finally, we have projection matching methods \cite{proj1,proj2,proj3} that determine the projection geometry and reconstruction jointly. Projection matching methods work by framing the CT reconstruction problem as a joint nonlinear least-squares optimization problem for the unknown attenuation \textit{and} the projection geometry. Due to the nonlinear coupling of the reconstruction and the projection geometry the joint optimization problem is often solved using an alternating approach, where either the reconstruction or projection geometry is fixed, while solving for the other quantity. This alternating approach is framed into an outer iterative scheme, where each iteration updates the reconstruction and projection geometry. However there is no guarantee for convergence using this alternating scheme, and many reconstructions may need to be computed to accurately estimate the projection geometry. 
\newline
\newline
A common theme for the described methods is that they only yield a point-estimate of the corrected projection geometry --- they do not provide a statistical measure of how (un)certain the computed estimate is. One contribution of our method is to address this limitation. Our proposed method is based on a Bayesian framework, where we determine the projection geometry and the reconstruction jointly by approximately sampling from a posterior probability distribution using a hierarchical Metropolis-within-Gibbs sampler. In this way, we naturally obtain uncertainty estimates of the projection geometry. Through numerical experiments with real tomographic data, we illustrate that the performance of the proposed method is similar to two existing geometry correcting methods for high-quality low noise data, and that we have superior performance under challenging conditions with fewer projections or highly noisy data.
\newline
\newline
In this paper, we focus on correcting the center-of-rotation offset in 2D fan-beam CT, since this setup is common for many other geometry correcting methods studied in the literature. However, we emphasize that the proposed method is flexible by nature, and that it can be generalized to other CT setups, including parallel-beam and 3D CT setups.

\section{Data description}
\subsection{SparseBeads data}
We use real tomographic data from the SparseBeads collection \cite{Beads}. The data was collected using a 320/225 kV Nikon XTEK Bay micro-CT instrument, located at the Henry Mosely X-ray Imaging Facility, The University of Manchester. The scanned objects are plastic cylinders containing glass beads. The detector has $2000\times2000$ pixels, and the size of a single detector pixel is $0.2 \text{mm}\times 0.2 \text{mm}$. The lab CT system employed a cone-beam setup with 2520 equidistant projections acquired over $360\degree$, and the source-to-center distance and center-to-detector distance were measured as $122$mm and $1400$mm respectively. A single central slice was provided for each data set, corresponding to a fan-beam setup with sinograms of size $2520\times 2000$. Using terminology from \cite{Beads}, we consider the $\mathrm{B3L1}$ and $\mathrm{Low\_B3L1}$ data sets in this paper. These data sets are scans of the same object, but with different exposure time per projection, $1000$ms for $\mathrm{B3L1}$ and $134$ms for $\mathrm{Low\_B3L1}$, so we have access to a high-dose and low-dose data set of the same object.

\section{Extended Reconstruction Model}
\subsection{Center-of-rotation misalignment for fan-beam CT}
\begin{figure}[ht]\centering
\includegraphics[width=0.50\textwidth]{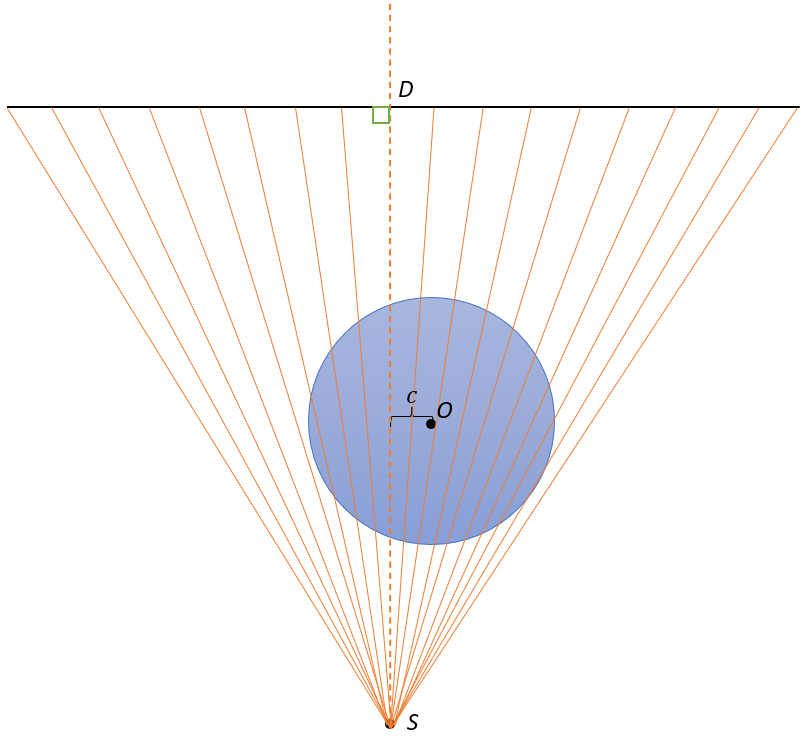}\hspace{0.05\textwidth}%
\begin{minipage}[b]{0.45\textwidth}\caption{Illustration of geometric misalignment in fan-beam CT. The orange lines illustrate the radiation emitting from the source $S$ penetrating the blue object domain with center-of-rotation $O$, and the horizontal black line illustrates the detector. The midline is the dotted line, which intersects with the detector perpendicularly at the point $D$. The center-of-rotation offset, $c$, is defined as the signed distance from the midline to the center-of-rotation $O$. The illustration depicts a configuration where $c$ is positive.}
\label{fig:misalignment}
\end{minipage}
\end{figure}

We consider fan-beam CT setups with a linear detector. We regard the horizontal offset of the center-of-rotation from the midline as an unknown quantity, see figure \ref{fig:misalignment}. We parametrize the center-of-rotation offset with the parameter $c$, which is defined as the signed distance from the midline to the center-of-rotation, and we define the sign such that figure \ref{fig:misalignment} illustrates a situation where $c>0$. For $c = 0$ the midline exactly penetrates the center-of-rotation.

\begin{figure}
    \centering
    \includegraphics[width = 1\textwidth]{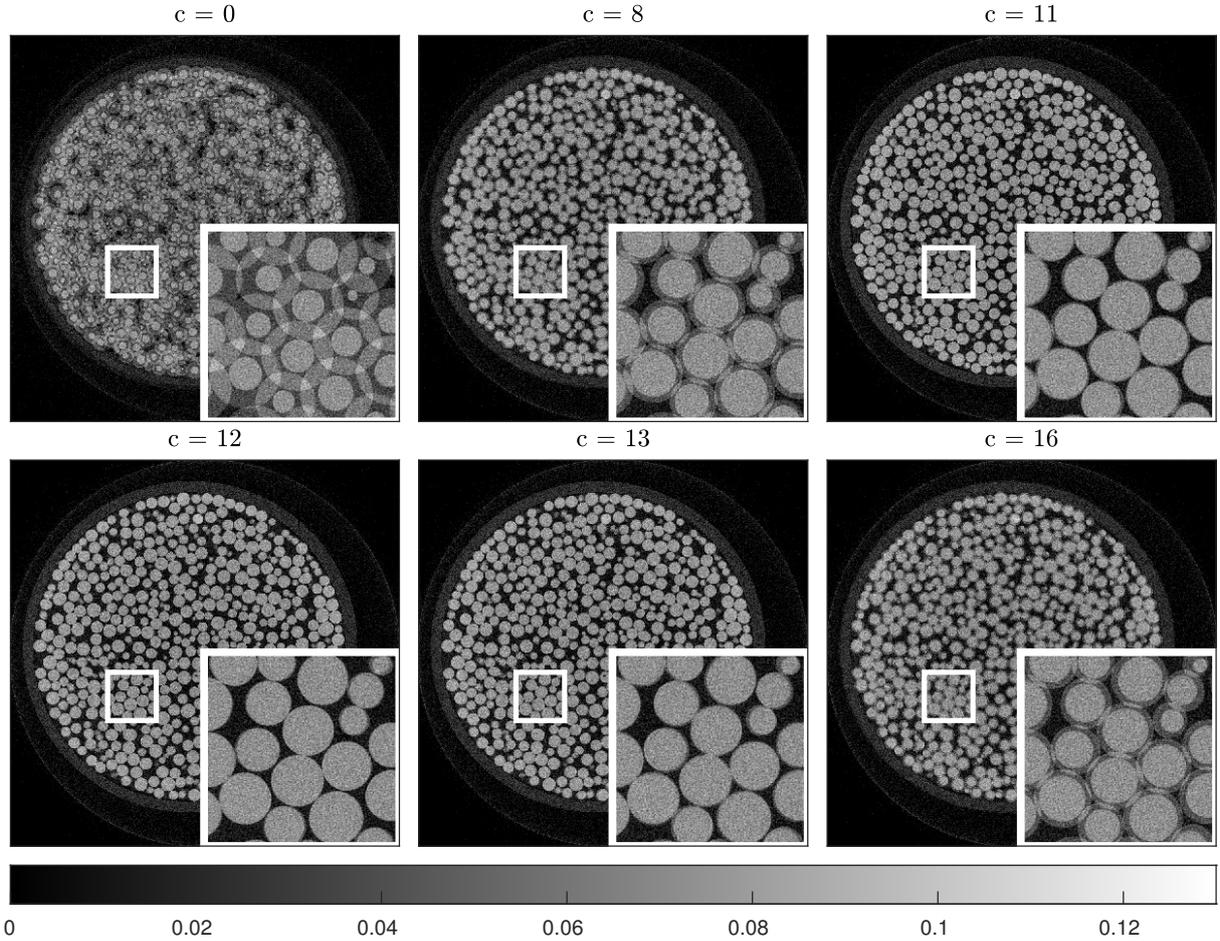}
    \caption{Illustration of reconstructions for high-dose data using different values of the center-of-rotation offset. Based on visual inspection the optimal value is close to $c = 12$.}
    \label{fig:BeadsRecon}
\end{figure}

In an idealized fan-beam setup we have $c = 0$, but in practice the center-of-rotation may be significantly offset from the midline. A failure to accurately estimate the center-of-rotation offset in the reconstruction model may lead to reconstructions with considerable artifacts. We illustrate the effect that the center-of-rotation offset has on the reconstruction quality in figure \ref{fig:BeadsRecon}. The $2000 \times 2000$ pixel reconstructions (pixel side length $17.4\si{\micro\metre}$) were computed using the high-dose \text{B3L1} data, and we measure $c$ in units of pixels such that $c = 1$ is a center-of-rotation offset of $17.4\si{\micro\metre}$ in physical units. For a given value of $c$, we frame the reconstruction problem as a Tikhonov regularized least-squares optimization problem with nonnegativity constraints
\begin{equation}
\label{eq:optimproblem}
    \underset{x\geq 0}{\text{min}}\quad ||A_c x - b||_2^2 + \alpha||x||_2^2,
\end{equation}
where $b\in \mathbb{R}^m$ is the measured sinogram, $x\in \mathbb{R}^n$ is the absorption coefficients of the object, $A_{c}\in \mathbb{R}^{m\times n}$ is the linear forward projection operator parametrized by the center-of-rotation offset $c\in \mathbb{R}$, and $\alpha\geq 0$ is the regularization parameter. We solve \eqref{eq:optimproblem} using the Fast Iterative Shrinkage-Thresholding Algorithm (FISTA), where the smooth part is the sum of the two terms in \eqref{eq:optimproblem} and the non-smooth part is the indicator function for the nonnegative orthant, see \cite{fista} for details. We use the ASTRA toolbox \cite{astra1,astra2} for the forward and back projections needed for FISTA. As we are dealing with a real data set, we do not know the true value of the center-of-rotation offset, but based on visual inspection of the reconstructions in figure \ref{fig:BeadsRecon} the true center-of-rotation offset is close to $12$ pixels for this particular data set, and we observe that relatively small deviations from this value result in considerable reconstruction artifacts that reduce the quality of the reconstruction.

\subsection{Bayesian model}
We consider an extended reconstruction model, which explicitly takes uncertainty in the center-of-rotation offset into account
\begin{equation}
    \label{eq:ExtendedModel}
    b = A_{c} x + \epsilon,
\end{equation}
where $\epsilon\in\mathbb{R}^m$ is additive noise and the remaining quantities are as defined in \eqref{eq:optimproblem}. Throughout we assume that the noise is Gaussian $\epsilon \sim N(0,\lambda^{-1} I)$, where $\lambda>0$ is the precision (inverse variance) of the noise. Looking at problem \eqref{eq:ExtendedModel} from a Bayesian perspective, we consider $x$ and $c$ as independent random variables, so that the likelihood of \eqref{eq:ExtendedModel} is given by 
\begin{equation}
\label{eq:likelihood}
    \pi(b|x,c) \propto \exp\left(-\frac{\lambda}{2} ||A_c x - b||_2^2\right).
\end{equation}
The maximum likelihood (ML) solution can be obtained by maximizing \eqref{eq:likelihood} over the variables $x$ and $c$. However, since CT is an ill-posed inverse problem, the ML solution will be dominated by high-frequency noise. To obtain high-quality reconstructions we must impose regularization, which in the Bayesian setting is done by considering prior distributions for the reconstruction and the projection geometry, that reflect our prior knowledge of the solution. Examples of common prior distributions in computational imaging include Gaussian \cite{Gaussian-MCMC}, Cauchy \cite{Cauchy-MCMC}, and edge-preserving priors such as total variation \cite{TV-MCMC}. The cost of using more informative priors is that increasingly sophisticated (inversion) methods must be used, which may lead to increased computational cost. The model \eqref{eq:ExtendedModel} has the added complexity of a stochastic forward operator, so we choose Gaussian priors $x|\delta \sim N(0,\delta^{-1} \Sigma_x)$ and $c \sim N(\mu_{c}, \sigma^2_{c})$ for simplicity. Here $\delta$ controls the intensity of the prior structure imposed by the precision matrix $\Sigma_x^{-1}$. The quantities $\delta$ and $\lambda$ (refered to as hyperparameters going forward) serve the role in the Bayesian inversion that the classical regularization parameter does in the non-Bayesian setting. We consider the hyperparameters as stochastic with conjugate Gamma hyperpriors $\lambda \sim \text{Gamma}(\alpha_{\lambda},\beta_{\lambda})$ and $\delta \sim \text{Gamma}(\alpha_{\delta},\beta_{\delta})$, where we assume that the shape parameters, $\alpha_{\delta}, \alpha_{\lambda}>0$, and inverse scale parameters, $\beta_{\delta}, \beta_{\lambda}>0$, are known. Regarding $\delta$ and $\lambda$ as stochastic hyperparameters imposes a natural hierarchy in the model and for this reason models with this structure are often denoted as hierarchical Bayesian models, see e.g. \cite{BardsleyBook}.
\newline
\newline
In the Bayesian setting, the solution to the inverse problem \eqref{eq:ExtendedModel} is a high-dimensional posterior probability distribution $\pi(x,c,\lambda,\delta|b)$. Using Bayes' law and assuming independence of the variables, the posterior distribution is given by
\begin{equation*}
    \pi(x,c,\lambda,\delta|b)\propto \pi(b|x,c,\lambda) \pi(x|\delta)\pi(c)\pi(\delta)\pi(\lambda),
\end{equation*}
which for our concrete likelihood and priors become
\begin{align}
\label{eq:posterior}
    \begin{split}
    \pi(x,c,\lambda,\delta|b)\propto\ &\lambda^{m/2 + \alpha_{\lambda} - 1} \delta^{n/2 + \alpha_{\delta} - 1} \\ &\exp\left(-\lambda\left(\frac{1}{2}||A_{c}x-b||_2^2+\beta_{\lambda}\right)-\delta\left(\frac{1}{2}||x||^2_{\Sigma_{x}^{-1}}+\beta_{\delta}\right)-\frac{1}{2\sigma_c^2}(c-\mu_{c})^2\right),
    \end{split}
\end{align}
where $||x||^2_{\Sigma_x^{-1}} = x^T \Sigma_x^{-1} x$. The benefits of using a Gaussian prior for the reconstruction is simplicity and computational efficiency, since we have a closed-form expression for the posterior distribution. However, it is also a highly uninformative prior that does not reflect much of our prior knowledge of the solution. For this reason, we also consider an augmentation to the Gaussian prior, where we add an implicit nonnegativity assumption for the reconstruction. Nonnegativity is a natural prior in CT, since the absorption coefficients are nonnegative by nature, and taking this into account in the inversion often results in higher quality reconstructions. To incorporate nonnegativity in the prior, we consider a mixture between two distributions: a truncated Gaussian defined on the positive orthant, and a distribution on the boundary of the positive orthant, which is obtained by appropriately transporting probability mass from outside the positive orthant. The resulting posterior distribution using this prior no longer has a simple closed-form expression, however \cite{nonneg1,nonneg2} show that posterior distribution samples using this prior can be obtained by solving a constrained least-squares optimization problem, see next section for more details.

\section{Hierarchical Metropolis-within-Gibbs sampler}
One way to use the posterior distribution to solve the inverse problem is to compute a maximum a posteriori (MAP) reconstruction by solving a nonlinear optimization problem. However, this does not provide us with any estimates of the uncertainty. Furthermore, the fact that the forward projection operator depends explicitly on the center-of-rotation offset $c$ complicates the resulting optimization problem significantly, since $A_c$ is a highly nonlinear function of $c$. For these reasons, we take a sampling approach, and we will develop a hierarchical Metropolis-within-Gibbs sampler to approximately sample from the full posterior distribution. A Gibbs sampler samples from the posterior distribution by iteratively sampling from all the conditional distributions. The overall structure of the sampler we propose is inspired from the approach presented in \cite{BardsleyBook}.

\subsection{Conditional distributions}
To efficiently use a Gibbs sampler to sample from the posterior distribution, we need efficient ways of sampling from the conditionals of the posterior distribution (denoted henceforth as conditional posterior distributions). The conditional posterior distributions for the hyperparameters are
\begin{subequations}
\begin{align}
    \pi(\lambda|x,b,c)&\propto \lambda^{m/2 + \alpha_{\lambda} - 1}\exp\left(-\lambda\left(\frac{1}{2}||A_cx - b||_2^2 + \beta_{\lambda}\right)\right), \label{eq:lambdacondpost}\\
    \pi(\delta|x,b)&\propto \delta^{n/2 + \alpha_{\delta} - 1}\exp\left(-\delta\left(\frac{1}{2}||x||^2_{\Sigma_x^{-1}} + \beta_{\delta}\right)\right). \label{eq:deltacondpost}
\end{align}
\end{subequations}
We immediately recognize these as Gamma distributions
\begin{subequations}
\begin{align*}
    \lambda|x,b,c &\sim \text{Gamma}\left(m/2 + \alpha_{\lambda},\frac{1}{2}||A_{c} x - b||_2^2 + \beta_{\lambda}\right),\\
    \delta|x,b&\sim \text{Gamma} \left(n/2 + \alpha_{\delta},\frac{1}{2}||x||_{\Sigma_x^{-1}}^2 + \beta_{\delta}\right).
\end{align*}
\end{subequations}
The conditional posterior distribution for the reconstruction is given by
\begin{equation*}
    \pi(x|\lambda,\delta,b,c)\propto \exp\left(-\frac{1}{2}\left(\lambda||A_{c}x - b||_2^2 + \delta ||x||_{\Sigma_x^{-1}}^2\right)\right).
\end{equation*}
Algebraic manipulations show that this is multivariate Gaussian
\begin{equation}
\label{eq:x-conditional}
    x|\lambda,\delta,b,c\sim N\left(\left(\lambda A_{c}^T A_{c} + \delta \Sigma_x^{-1}\right)^{-1}\lambda A_{c}^{T} b,\left(\lambda A^T_{c}A_{c} + \delta \Sigma_x^{-1}\right)^{-1}\right).
\end{equation}
Sampling from a high-dimensional multivariate Gaussian distribution is generally expensive, however the distribution coincides with the family of solutions to a stochastic regularized least-squares problem
\begin{equation}
\label{eq:Gaussianopt}
    x(\xi_m,\xi_n) = \underset{x}{\text{argmin}} \quad \left\{\frac{\lambda}{2}||A_{c}x-b-\lambda^{-1/2} \xi_m||_2^2 + \frac{\delta}{2}||\Sigma_{x}^{-1/2}x-\delta^{-1/2} \xi_n||_2^2\right\},
\end{equation}
where $\xi_m \in \mathbb{R}^m$ and $\xi_n \in \mathbb{R}^n$ are Gaussian \textit{i.i.d.} samples, see \cite{BardsleyBook} for details. This suggests that approximate samples from the multivariate Gaussian distribution can be obtained by inexact solutions of the stochastic optimization problem.
\newline
\newline
If we instead want to obtain samples using the implicit nonnegativity prior, we follow \cite{nonneg1,nonneg2} and define the conditional posterior for the reconstruction as the family of solutions to the constrained stochastic regularized least-squares problem
\begin{equation}
\label{eq:Gaussiannonnegopt}
    x(\xi_m,\xi_n) = \underset{x\geq 0}{\text{argmin}} \quad \left\{\frac{\lambda}{2}||A_{c}x-b-\lambda^{-1/2} \xi_m||_2^2 + \frac{\delta}{2}||\Sigma_{x}^{-1/2}x-\delta^{-1/2} \xi_n||_2^2\right\}.
\end{equation}
The only difference compared to the model without nonnegativity is the addition of the constraint. Again following \cite{nonneg1,nonneg2}, the conditional prior for the hyperparameter $\delta$ when using the implicit nonnegativity prior takes the form
\begin{equation}
\label{eq:deltanonneg}
    \pi(\delta|x,b) = \delta^{\bar{n}/2 + \alpha_{\delta}-1} \exp\left(-\delta\left(\frac{1}{2}||x||_{\Sigma_x^{-1}} + \beta_x\right)\right),
\end{equation}
where $\bar{n}$ is the number of non-zero entries in $x$.
\newline
\newline
Finally the conditional posterior distribution for $c$ is
\begin{equation}
\label{eq:gemetriccond}
    \pi(c|x,\lambda,b) \propto \exp\left(-\frac{\lambda}{2}||A_{c}x - b||_2^2 - \frac{1}{2 \sigma_c^2}(c - \mu_{c})^2\right).
\end{equation}
This is not a simple closed-form distribution due to the stochastic forward operator $A_{c}$, so we must use a Markov Chain Monte Carlo (MCMC) method to obtain samples from this conditional distribution.

\subsection{Implementation}
To obtain a sample from the conditional posterior of the reconstruction, we must solve \eqref{eq:Gaussianopt} (or \eqref{eq:Gaussiannonnegopt} in the case of nonnegativity constraints) with a high level of accuracy. However, this is too costly for high-dimensional problems, where we may need thousands of samples. To address this issue and increase computational efficiency, we instead use FISTA with a (small) fixed number of iterations to obtain approximate solutions to \eqref{eq:Gaussianopt} (or \eqref{eq:Gaussiannonnegopt} in case of nonnegativity constraints). Early termination of FISTA incurs an approximation error in the conditional sample of the reconstruction, which implies that we can not guarantee stationarity of the full Metropolis-within-Gibbs sampler. To mitigate some of the error caused by truncated iterations, we use the previous sample as the initial guess for the FISTA iterations as a warm-start procedure. Truncation of an iterative method to sample from a multivariate Gaussian resembles the Gradient Scan method presented in \cite{GradientScan}. The difference is that the Gradient Scan method uses a fixed number of conjugate gradient iterations instead of FISTA iterations. Using properties of the conjugate gradient iterates the authors of \cite{GradientScan} prove asymptotic convergence. However, their proof does not extend to the case with nonnegativity constraints. Using FISTA instead of conjugate gradient has the benefit of being easily adapted to solve the nonnegativity case \eqref{eq:Gaussiannonnegopt}.
\newline
\newline
To obtain samples from the conditional posterior distribution for the center-of-rotation offset, we use Metropolis-Hastings, see e.g. \cite{MarkovSpringer}, with a symmetric Gaussian random walk proposal distribution
\begin{equation*}
    q(\hat{c}|c) \propto \exp\left(-\frac{1}{2s^2}(\hat{c} - c)^2\right),
\end{equation*}
where $s>0$ is the \textit{step size} in the Metropolis-Hastings method, which can be tuned to obtain a reasonable acceptance rate in the Metropolis-Hastings sampling procedure. The optimal acceptance rate is highly dependent on the problem, but we generally aim for around $25\%$ as is also discussed in \cite{MetroAcceptance}. The Metropolis-Hastings algorithm is shown in Algorithm \ref{alg:MH}.

\begin{algorithm}
    \caption{Metropolis-Hastings for center-of-rotation offset} \label{alg:MH}
    \begin{algorithmic}
        \STATE  \textbf{Input}: reconstruction, $x\in \mathbb{R}^{n}$; initial value of center-of-rotation offset, $c^0\in \mathbb{R}$; center-of-rotation offset Gaussian mean, $\mu_c\in \mathbb{R}$; center-of-rotation offset Gaussian standard deviation, $\sigma_c>0$; Metropolis-Hastings step size, $s>0$; number of Metropolis-Hastings iterations, $k_{\mathrm{metro}}$; noise hyperparameter, $\lambda>0$; sinogram, $b\in \mathbb{R}^{m}$.
        \STATE \textbf{Output}: Sample from $\pi(c|x,\lambda,b)$.
        \FOR{$k=1$ to $k_{\mathrm{metro}}$}
        \STATE Sample from Gaussian proposal: $c^{\mathrm{prop}} \sim N(c^{k-1},s^2)$
        \STATE Compute acceptance probability: $p^{\text{accept}} = \text{min}\left\{1,\frac{\pi(c^{\mathrm{prop}}|x,\lambda,b)}{\pi(c^{k-1}|x,\lambda,b)}\right\}$
        \STATE Set $c^k = c^{\mathrm{prop}}$ with probability $p^{\text{accept}}$. Otherwise let $c^k = c^{k-1}$.
        \ENDFOR
    \end{algorithmic}
\end{algorithm}
As we use Metropolis-Hastings to sample from one of the conditional distributions our proposed sampler for the full posterior distribution is a hierarchical Metropolis-within-Gibbs sampler, see e.g. \cite{MetroWithinGibbs,MarkovReview} for details on Metropolis-within-Gibbs samplers. The full sampler using the prior without nonnegativity is presented in Algorithm \ref{alg:HierarchicalMetroWithinGibbs}.

\begin{algorithm}[ht]
    \caption{Hierarchical Metropolis-within-Gibbs Sampler}
    \label{alg:HierarchicalMetroWithinGibbs}
    \begin{algorithmic}
        \STATE \textbf{Input}: initial reconstruction, $x^0\in \mathbb{R}^n$; initial value of center-of-rotation offset, $c^0\in \mathbb{R}$; center-of-rotation offset Gaussian mean, $\mu_{c} \in \mathbb{R}$; center-of-rotation offset Gaussian standard deviation, $\sigma_{c}>0$; Metropolis-Hastings step size, $s>0$; number of Gibbs samples, $k_{\mathrm{Gibbs}}$; number of Metropolis-Hastings iterations, $k_{\mathrm{Metro}}$;  number of FISTA iterations, $k_{\mathrm{FISTA}}$; Gamma distribution parameters, $\gamma_{\lambda}$, $\gamma_{\delta}$, $\beta_{\lambda}$, $\beta_{\delta}>0$; sinogram, $b\in \mathbb{R}^m$.
        \STATE \textbf{Output}: approximate sample chain $\{(x^k,c^k,\lambda^k,\delta^k)\}_{k=1}^{k_{\mathrm{Gibbs}}}$ from full posterior distribution $\pi(x,c,\lambda,\delta|b)$.
        \FOR{$k=1$ to $k_{\mathrm{Gibbs}}$}
        \STATE \textit{Sample conditional hyperparameters from conjugate Gamma distributions}
        \STATE $\lambda^k|x^{k-1},c^{k-1},b \sim \text{Gamma}\left(m/2 + \alpha_{\lambda},1/2||A_{c^{k-1}} x^{k-1} - b||_2^2 + \beta_{\lambda}\right)$
        \STATE $\delta^k|x^{k-1},b\sim \text{Gamma} \left(n/2 + \alpha_{\delta},1/2||x||_{\Sigma_x^{-1}}^2 + \beta_{\delta}\right)$
        \STATE
        \STATE \textit{Sample conditional center-of-rotation offset using Metropolis-Hastings}
        \STATE Sample $c^k|x^{k-1},\lambda^k,\delta^k,b$ using $k_{\mathrm{Metro}}$ iterations of Algorithm 1
        \STATE
        \STATE \textit{Sample conditional reconstruction using truncated FISTA}
        \STATE Approximate $x^k|c^k,\lambda^k,\delta^k,b$ by running $k_{\mathrm{FISTA}}$ FISTA iterations on \eqref{eq:Gaussianopt} with initial condition $x^{[0]} = x^{k-1}$.
        \ENDFOR
    \end{algorithmic}
\end{algorithm}

Using the nonnegativity prior, Algorithm \ref{alg:HierarchicalMetroWithinGibbs} must be changed in the following way: the conditional $\delta$-sample step must be changed to \eqref{eq:deltanonneg}, and the conditional reconstruction sample steps must instead run on \eqref{eq:Gaussiannonnegopt}. The remaining steps are unchanged.
\newline
\newline
A good choice for the initial reconstruction $x^0\in \mathbb{R}^n$ can be obtained by solving a (nonnegativity-constrained) regularized least-squares problem like \eqref{eq:optimproblem} using the initial center-of-rotation offset $c^0$. In theory we only need to take a single Metropolis-Hastings iteration, $k_{\mathrm{Metro}} = 1$, for each full Gibbs iteration \cite{MetroWithinGibbs1}. However, we found that taking multiple Metropolis-Hastings steps speeds up the estimation of the center-of-rotation offset considerably. The reason for this is that taking multiple Metropolis-Hastings steps allows for potentially larger steps for each full Gibbs iteration, while also reducing the risk of complete rejections when sampling the center-of-rotation offset. Complete rejections are undesirable, since it is considerably more expensive to obtain a sample of the reconstruction compared to a sample of the center-of-rotation offset. In practice we found that $k_{\mathrm{Metro}} = 10$ yielded a nice balance between avoiding rejections, and also not using too many computational resources on the Metropolis-Hastings step. For the number of FISTA iterations, we use $k_{\mathrm{FISTA}} = 20$ unless otherwise stated, since we empirically found this to provide a reasonable balance between performance and computational effort.
\newline
\newline
\textbf{Computational remarks}. The computationally demanding steps of Algorithm \ref{alg:HierarchicalMetroWithinGibbs} is sampling the center-of-rotation offset and sampling the reconstruction. Sampling the center-of-rotation offset requires $k_{\mathrm{Metro}}$ forward projections, and sampling the reconstruction requires $k_{\mathrm{FISTA}}$ forward and back projections. The cost of a forward projection and a back projection is approximately equal, so the total cost of the algorithm is $k_{\mathrm{Gibbs}}(2k_{\mathrm{FISTA}} + k_{\mathrm{Metro}})$ forward projections. With our specific values for $k_{\mathrm{FISTA}}$ and $k_{\mathrm{Metro}}$, the total cost of the MCMC algorithm is $2.5k_{\mathrm{Gibbs}}k_{\mathrm{FISTA}}$ forward projections. Hence the computational scaling of the algorithm is similar to computing a MAP estimate using an iterative method, although the constant $2.5 k_{\mathrm{Gibbs}}$ may need to be large to obtain a satisfying reconstruction with a good estimate of the center-of-rotation offset. We remark that the proposed method also samples the regularization parameter, which implies that we do not need to use additional resources to tune this quantity, which must be done for classical MAP reconstruction.  

\section{Numerical Experiments}

In this section, we present numerical experiments using the proposed MCMC method. All experiments shown were conducted using the model with the implicit nonnegativity assumption, since this yielded superior results compared to the model without nonnegativity. We compare the MCMC method with two other commonly used methods, both of which are computationally cheap and easy to implement, since they are based solely on the recorded sinogram. The first method, COM, is from \cite{CenterOfMass}, and it is based on a center-of-mass calculation of the sinogram. The second method, XCORR, is from \cite{CrossCor}, and it is based on auto-correlation of the sum of projections in the sinogram. Both the COM and XCORR methods make a correction based on symmetry features in the sinogram. For this reason, they both require $360 \degree$ projection data, whereas the MCMC method makes no assumption regarding angular range of the projections.
\newline
\newline
We conduct numerical experiments for both the full $360\degree$ high-dose and low-dose data. For the high-dose data, we also conduct numerical experiments using a fast-scan setup, where we limit the number of projections, so that we do not have full $360\degree$ projection data, to illustrate the versatility of the proposed method. The results consist of reconstructions as well as sample chains for the MCMC method to see convergence of the hyperparameters and center-of-rotation offset. Overall the results show that the MCMC method performs similarly to the COM and XCORR method for the high-dose $360\degree$ data, and that we have superior performance for the more challenging low-dose and fast-scan setups. 

\subsection{High-dose data}
\begin{figure}
    \centering
    \includegraphics[width=\textwidth]{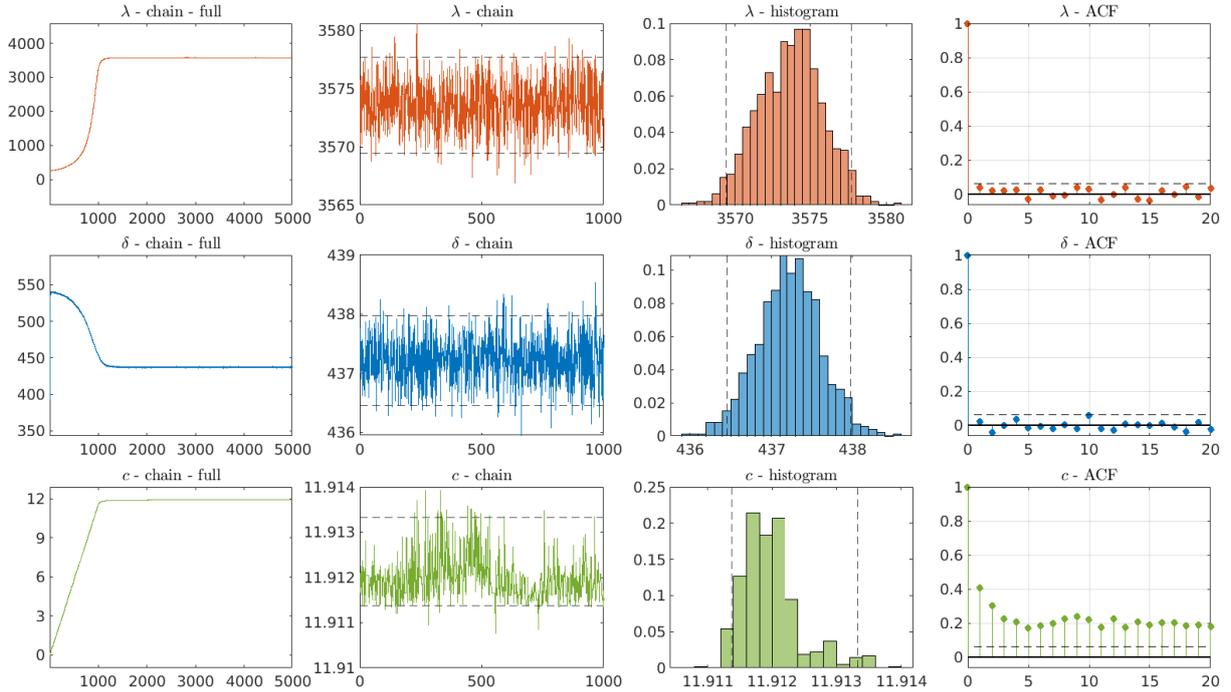}
    \caption{MCMC chains for high-dose SparseBeads data. The first column show the full MCMC chains, the second column show the chains after discarding the first $4000$ samples as burn-in, the third column show histograms after burn-in, and the final column is the autocorrelation function of the samples after burn-in. The dotted black lines indicate empirical $95\%$ credibility intervals. Based on visual inspection all chains look stationary after burn-in.}
    \label{fig:highdosebeads-chain}
\end{figure}

We use $c^0=0$, $s=2.5\cdot10^{-3}$, $\mu_c=0$, $\sigma_c=20$ (all measured in pixels), and we use the identity matrix for the Gaussian precision matrix $\Sigma_x^{-1}$. For the hyperprior parameters we choose $\alpha_{\lambda}=\alpha_{\delta} = 1$ and $\beta_{\lambda}=\beta_{\delta} = 10^{-4}$ as recommended in \cite{BardsleyBook}, yielding weakly informative exponential hyperpriors. We take a total of $k_{\mathrm{Gibbs}} = 5000$ samples, and we do not do any thinning. The resulting MCMC-chains can be seen in figure \ref{fig:highdosebeads-chain}; the full chains have three distinct phases - a fast transient phase for the first roughly $1500$ samples, where the hyperparameters and center-of-rotation offset rapidly change, a slower transient phase for the next $2500$ samples, where the samples slowly approach the equilibrium, and finally an equilibrium phase for the final $1000$ samples, where we obtain approximate samples from the posterior distribution. Based on visual inspection, the chains look stationary after discarding the first $4000$ samples for burn-in. We note that the autocorrelation for the center-of-rotation offset chain is higher than for the hyperparameters, which is due to the fact that we use Metropolis-Hastings to sample this. The center-of-rotation offset in the stationary phase is around 12 pixels, which aligns with the observation made in figure \ref{fig:BeadsRecon}. In addition, we observe that the variance of the $c$-chain is relatively low. Looking at the histogram the distribution of the center-of-rotation offset looks slightly skewed, but this is likely due to the relatively few samples taken, and we expect a bell-shaped distribution as we increase the number of samples.

The first column of figure \ref{fig:reconsallmethods} shows the reconstructions using the initial, the MCMC corrected, the COM corrected, and the XCORR corrected center-of-rotation offset for the high-dose data. The MCMC reconstruction is the posterior sample mean after the burn-in phase. The reconstructions for the initial, the COM, and the XCORR center-of-rotation offsets were obtained by solving \eqref{eq:optimproblem} to a high level of accuracy, where the regularization parameter was determined from the means of the hyperparameter chains after the burn-in phase. This is a fair comparison, since the solution to \eqref{eq:optimproblem} is a MAP estimate for the simpler reconstruction problem, where $c$, $\lambda$, and $\delta$ are assumed known, and you let $\alpha := \delta/\lambda$. We see a vast increase in reconstruction quality from all the correction methods compared to the initial estimate $c = 0$. The MCMC and XCORR reconstructions are slightly sharper compared to the COM method, but all methods produce high-quality reconstructions. 

\subsection{Low-dose data}
\begin{figure}
    \centering
    \includegraphics[width = \textwidth]{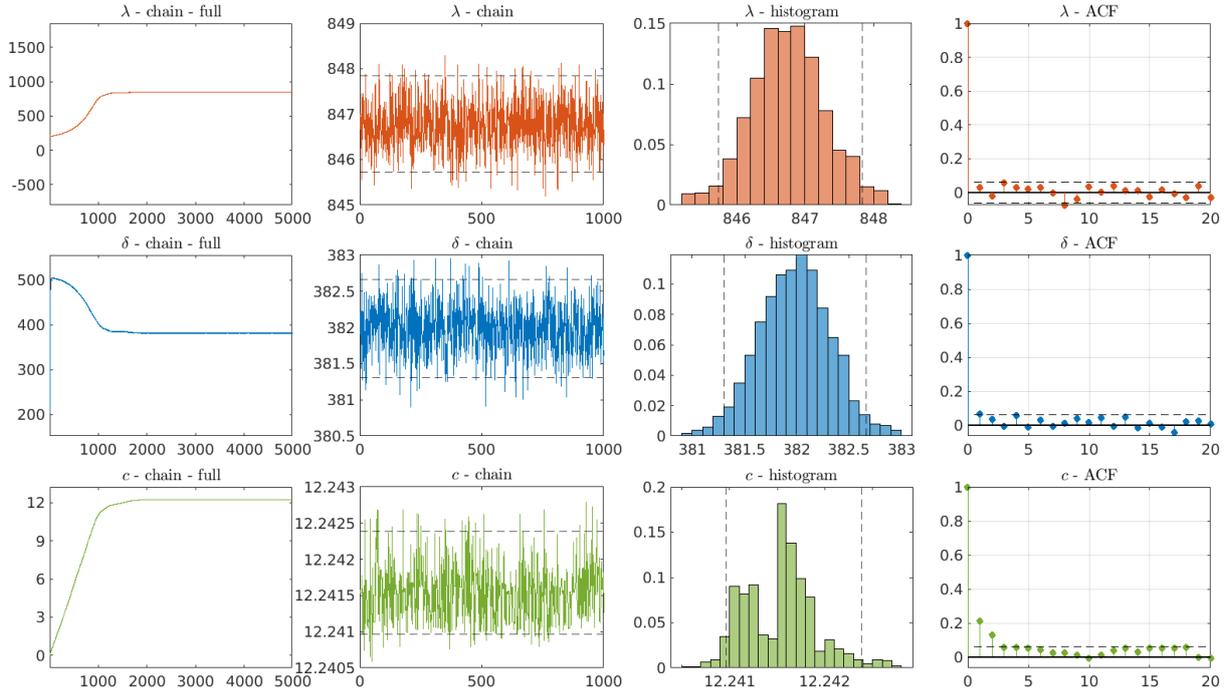}
    \caption{MCMC chains for low-dose SparseBeads data. Based on visual inspection the chains after burn-in look stationary. Notice that the values for the inverse noise variance $\lambda$ is much smaller compared to the high-dose case, reflecting the higher noise levels in the low-dose data.}
    \label{fig:lowdosebeads-chain}
\end{figure}

We use the same algorithm parameters as before, and the MCMC chains are shown in figure \ref{fig:lowdosebeads-chain}. The biggest difference between the high-dose and the low-dose MCMC chains is that the sampled $\lambda$-values are much smaller in the low-dose case. This is expected, since $\lambda$ represents the inverse variance of the noise, and we have more noise for the low-dose data. Interestingly, the noise level does not seem to affect the credibility interval of the $c$-estimate in this case, which is approximately the same width in the high-dose and low-dose case.

The second column in figure \ref{fig:reconsallmethods} shows the reconstructions for the low-dose data using the different geometry correction techniques. In general, we see that the reconstructions are noisier compared to the high-dose reconstructions, which is expected since the signal-to-noise ratio is lower for the low-dose data. All geometry correction methods demonstrate improved reconstruction quality compared to the initial center-of-rotation offset. However, contrary to the high-dose data, we see a significant improvement in clarity and sharpness for the MCMC method compared to the COM and XCORR method. This indicates that the center-of-rotation offset from the MCMC method is closer to the actual value compared to the estimates produced by the other methods, and that the MCMC method is more robust against noisy data.

\subsection{Fast-scan data}
\begin{figure}
    \centering
    \includegraphics[width = \textwidth]{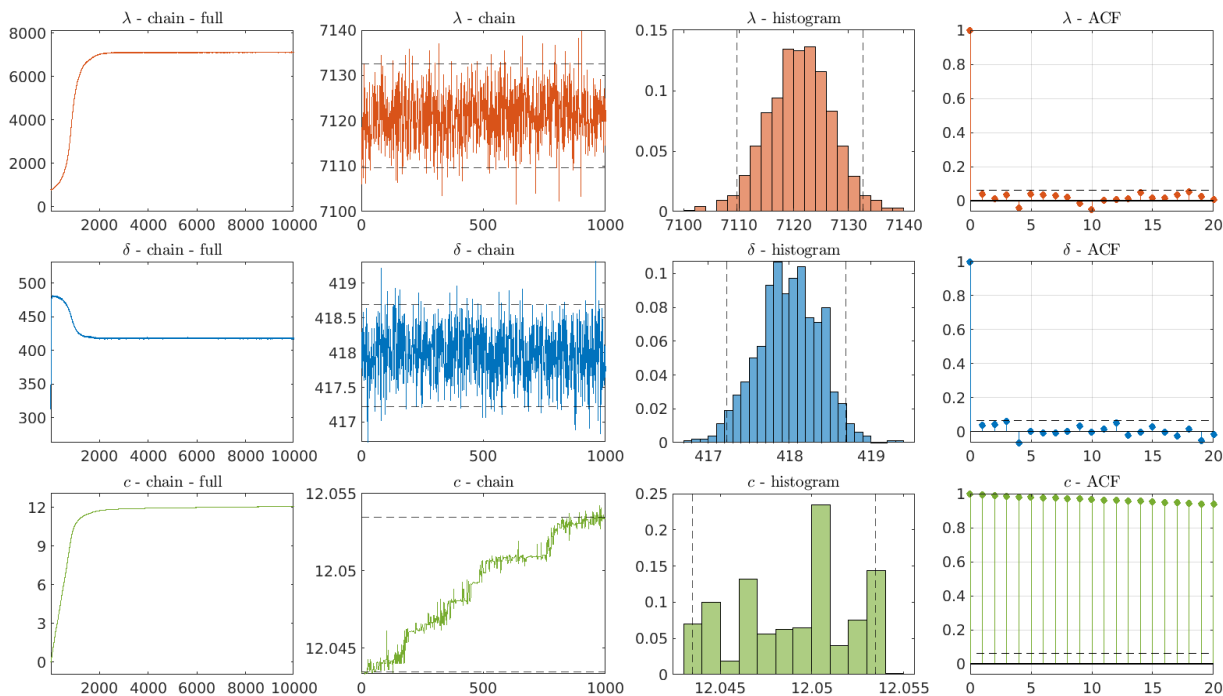}
    \caption{MCMC chains for high-dose fast-scan SparseBeads data. Notice that the $c$-chain does not look stationary in this case.}
    \label{fig:limitedanglechain}
\end{figure}

Finally, we conduct experiments, where we reduce the number of projections. We again use the high-dose data, but we restrict the projections to the interval $[0\degree,210\degree]$ to give a fast-scan setup. Thus the requirement regarding $360\degree$ projection data for the COM and XCORR method is not satisfied, so we can not expect good performance using these methods. We use the same algorithm parameters as for the previous experiments except that we take a total of $10000$ samples and discard the first $9000$. The resulting MCMC chain can be seen in figure \ref{fig:limitedanglechain}. We remark that the reduced angular coverage has an adverse effect on the apparent convergence rate of the MCMC chain, and the chain is still not stationary after $10000$ samples. However, we see very slow movement in the final $1000$ samples, and the estimate does look close to the optimal based on earlier experiments.

\begin{figure}
    \centering
    \includegraphics[width = 0.9\textwidth]{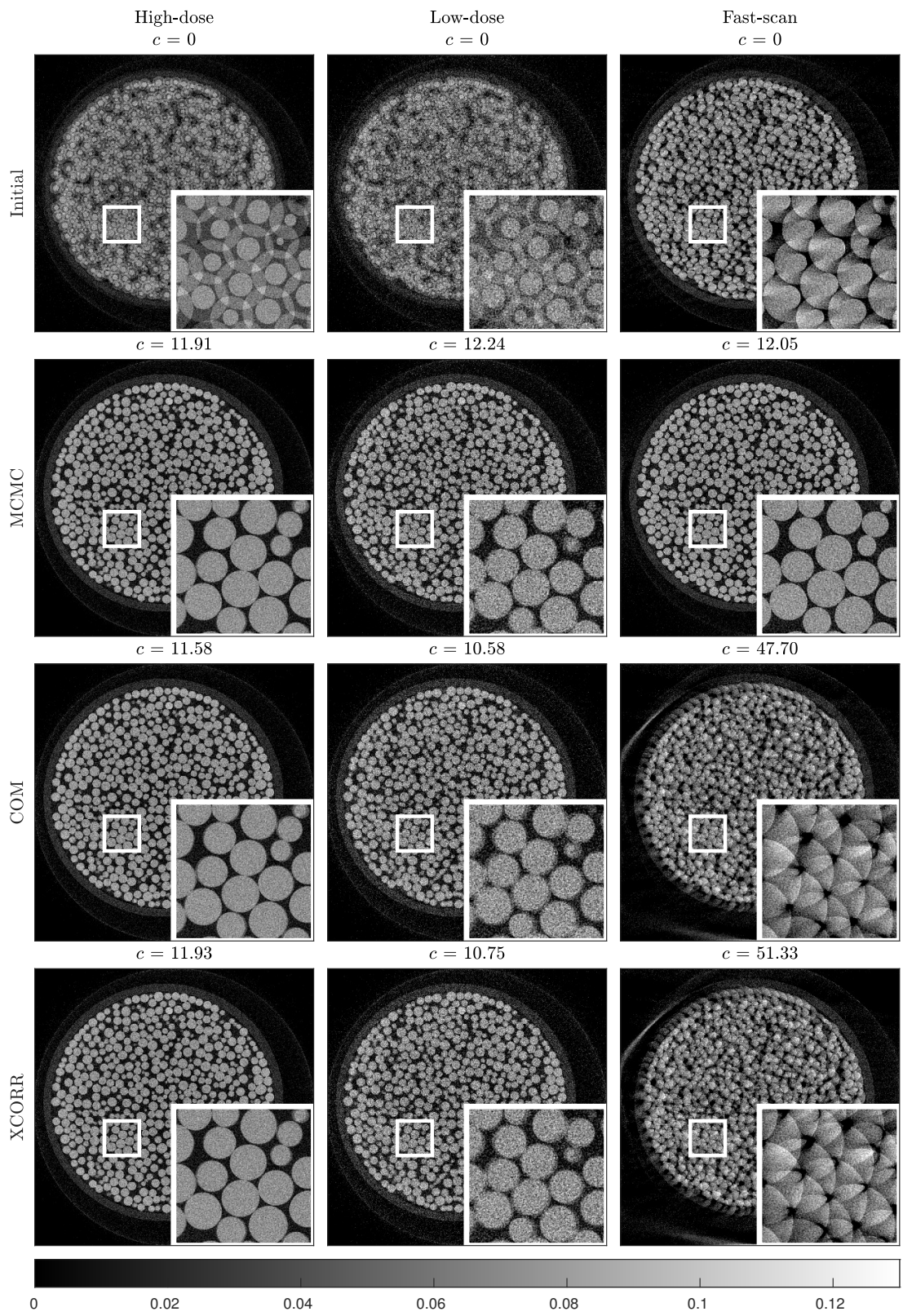}
    \caption{Reconstructions for different data sets and different geometry correction techniques.}
    \label{fig:reconsallmethods}
\end{figure}

The fast-scan reconstructions can be seen in the third column of figure \ref{fig:reconsallmethods}. Here we see a vast increase in quality of the reconstruction using the MCMC method compared to the COM and XCORR method, both of which fail to give a meaningful estimate. One may argue that the comparison in this case is not fair, since we do not satisfy the $360\degree$ requirement needed for the COM and XCORR method. However, we include the results to show challenging non-standard scenarios, where a more flexible and computationally expensive approach may be needed to obtain good results. Notice that the reconstruction for the MCMC method is slightly noisier compared to the full angle case due to the fewer projections in this setup.

\section{Discussion and Conclusion}
In CT reconstruction models, the projection geometry is often assumed to be known, whereas in reality it is uncertain and must be estimated. We have demonstrated that explicitly taking the uncertainty into account in an extended reconstruction model results in high-quality reconstructions under challenging conditions with real tomographic data. The focus in this paper has been on 2D fan-beam CT, but the proposed method is highly flexible, and it can easily be generalized to other geometries, for instance parallel-beam and 3D cone-beam CT, since we only need the ability to do forward and back projections with a given projection geometry. 
\newline
\newline
The biggest concern regarding the use of MCMC methods for high-dimensional inverse problems is the computational complexity. We have addressed this by developing an efficient Metropolis-within-Gibbs sampler that obtains approximate samples from the posterior distribution aided by modern optimization techniques. We stress that the proposed method is more expensive compared to many other geometry correcting methods, and for high-quality low noise data there is no need to use the expensive MCMC machinery, since the resulting estimates and reconstructions are comparable in quality to cheaper simpler methods. However, for challenging high-noise and/or fast-scan setups, the proposed method may yield superior results. Finally, we remark that simple geometry correction methods can work excellently in tandem with the proposed MCMC method --- we can obtain an excellent initial guess of the geometry using the simple methods, and then subsequently improve upon the estimate using the MCMC method. This may also significantly reduce the computational cost of the MCMC method, since we can skip a large part of the initial transient phase in the MCMC chains.
\newline
\newline
The proposed method can be improved in many ways. The main improvement is to extend the model to include more geometric parameters. In this paper we focus solely on the center-of-rotation offset for fan-beam CT. However, this single parameter does not completely determine the projection geometry for fan-beam setups. Instead, we should ideally fully control the horizontal and vertical position of the radiation source and detector. This would potentially allow us to correct more challenging misalignments in the projection geometry, which might result in higher quality reconstructions. The only change needed in the proposed algorithm is including more geometric variables in the Metropolis-Hastings step. 
\newline
\newline
The method can also be improved by considering more informative prior distributions instead of the Gaussian. This may further improve the quality of the reconstructions, and it may also assist in obtaining even more accurate estimates of the geometric parameters.

\section*{Acknowledgments}
This work was supported by The Novo Nordisk Foundation (Grant No. NNF20OC0061894) and by The Villum Foundation (Grant No. 25893).

\section*{References}

\bibliography{iopart-num}

\end{document}